\documentclass[12pt]{article}
\usepackage{graphicx, tikz, color,url}
\usepackage{amsmath,amssymb,latexsym}
\usetikzlibrary{arrows}
\usepackage{soul}
\usepackage{comment}

\newtheorem{theorem}{Theorem}[section]
\newtheorem{definition}[theorem]{Definition}
\newtheorem{lemma}[theorem]{Lemma}
\newtheorem{proposition}[theorem]{Proposition}

\newtheorem{question}[theorem]{Question}

\newtheorem{problem}[theorem]{Problem}

\newcommand{\proof}{\noindent{\bf Proof.\ }}
\newcommand{\qed}{\hfill $\square$ \bigskip}

\newcommand{\gsmb}{\gamma_{\rm SMB}}
\newcommand{\gmb}{\gamma_{\rm MB}}

\begin{document}
\title{Maker-Breaker domination game critical graphs}
\author{ Athira Divakaran$^{a,}$\thanks{Email: \texttt{athiradivakaran91@gmail.com}}
\and 
 Tanja Dravec$^{d,e,}$\thanks{Email: \texttt{tanja.dravec@um.si}}
\and
Tijo James$^{b,}$\thanks{Email: \texttt{tijojames@gmail.com}} 
\and 
Sandi Klav\v{z}ar$^{c,d,e,}$\thanks{Email: \texttt{sandi.klavzar@fmf.uni-lj.si}}
\and
Latha S Nair$^{a,}$\thanks{Email: \texttt{lathavichattu@gmail.com}} 
}
\maketitle

\begin{center}
$^a$ Department of Mathematics, Mar Athanasius College, Kothamangalam, India\\
	\medskip
$^b$ Department of Mathematics, Pavanatma College, Murickassery, India\\
\medskip
$^c$ Faculty of Mathematics and Physics, University of Ljubljana, Slovenia \\
\medskip

$^{d}$ Institute of Mathematics, Physics and Mechanics, Ljubljana, Slovenia \\
\medskip

$^e$ Faculty of Natural Sciences and Mathematics, University of Maribor, Slovenia \\
\medskip
\end{center}

\begin{abstract}
The Maker-Breaker domination game (MBD game) is a two-player game played on a graph $G$  by  Dominator and Staller. They alternately select unplayed vertices of $G$. The goal of Dominator is to form a dominating set with the set of vertices selected by him while that of Staller is to prevent this from happening. In this paper MBD game critical graphs are studied. Their existence is established and critical graphs are characterized for most of the cases in which the first player can win the game in one or two moves. 
\end{abstract}


\noindent
{\bf AMS Subj.\ Class.\ (2020)}: 05C57, 05C69

 \section{Introduction}
\label{sec: intro}

It is usually said that the Maker-Breaker game was introduced by Erd\H{o}s and Selfridge in~\cite{erdos-1973}, but it appeared already earlier in a paper of Hales and Jewett~\cite{HJ-1963}. Anyhow, this is a two-person game played on an arbitrary hypergraph $\mathcal{H}$. The players named Maker and Breaker alternately select an unplayed vertex of $\mathcal{H}$ during the game. Maker aims to occupy all the vertices of some hyperedge, on the other hand, Breaker's goal is to prevent Maker from doing it. The game has been extensively researched, both in general and in specific cases, cf.\ the book~\cite{hefetz-2014}, the recent paper~\cite{rahman-2023}, and references therein.

In this paper, we are interested in the domination version of the Maker-Breaker game which was introduced in 2020 by Duch\^{e}ne, Gledel, ~Parreau, and Renault~\cite{duchene-2020}. The {\em Maker-Breaker domination game} ({\em MBD game} for short) is a game played on a graph $G=(V(G), E(G))$  by two players named Dominator and Staller. These names were chosen so that the players are named in line with the previously intensively researched domination game~\cite{bresar-2010, book-2021}. Just as in the general case, the two players alternately select unplayed vertices of $G$. The aim of Dominator is to select all the vertices of some dominating set of $G$, while Staller aims to select at least one vertex from every dominating set of $G$. There are two variants of this game depending on which player has the first move. A {\em D-game} is the MBD game in which Dominator has the first move and an {\em S-game} is the MBD game in which Staller has the first move. 

The following graph invariants are naturally associated with the MBD game~\cite{bujtas-2024, gledel-2019}. The {\em Maker-Breaker domination number}, $\gmb(G)$, is the minimum number of moves of Dominator to win the D-game on $G$ when both players play optimally. That is, $\gmb(G)$ is the minimum number of moves of Dominator such that he wins in this number of moves no matter how Staller is playing. If Dominator has no winning strategy in the D-game, then set $\gmb(G)= \infty$. The {\em Staller-Maker-Breaker domination number},  $\gsmb(G)$, is the minimum number of moves of Staller to win the D-game on $G$ when both players play optimally, where $\gsmb(G)= \infty$ if Staller has no winning strategy.  In a similar manner, $\gmb'(G)$ and $\gsmb'(G)$ are the two parameters associated with the S-game. Briefly, we will call $\gmb(G)$, $\gmb'(G)$, $\gsmb(G)$, and $\gsmb'(G)$ the {\em MBD numbers of} $G$.

In the seminal paper~\cite{duchene-2020} it was proved, among other things, that deciding the winner of the MBD game can be solved efficiently on trees, but it is PSPACE-complete even for bipartite graphs and split graphs. In~\cite{duchene-2023+} the authors give a complex linear algorithm for the Maker-Maker version of the game played on forests. The paper~\cite{bujtas-2024} focuses on $\gsmb(G)$ and $\gsmb'(G)$ and among other results establishes an appealing exact formula for $\gsmb'(G)$ where $G$ is a path. In~\cite{bujtas-2023}, for every positive integer $k$, trees $T$ with $\gsmb(T) = k$ are characterized and exact formulas for $\gsmb(G)$ and $\gsmb'(G)$ derived for caterpillars. In the main result of~\cite{forcan-2023+}, $\gmb$ and $\gmb'$ are determined for Cartesian products of $K_2$ by a path. In~\cite{dokyeesun-2024+}, the MBD game is further studied on Cartesian products of paths, stars, and complete bipartite graphs. The total version of the MBD game was introduced in~\cite{gledel-2020} and further investigated in~\cite{forcan-2022}.

It is known from~\cite[Lemma~2.3]{dokyeesun-2024+}  that all the four graph invariants associated with the MBD game are monotonic for adding and deleting an edge. This motivates us to introduce MBD game critical graphs as respectively the graphs for which $\gmb$, $\gmb'$, $\gsmb$, $\gsmb'$ changes as soon as an arbitrary edge is removed/added. For the classical domination game, this aspect has already been studied  in~\cite{bujtas-2015, henning-2018, Kexiang-2018, wora-2024}. In the preliminaries,  additional definitions are listed and known results about the MBD game needed later on are stated. In Section~\ref{sec:main}, $\gmb$- and $\gmb'$-critical graphs are formally defined. In the same section the existence of the $\gmb$- and $\gmb'$-critical graphs is established and explicit constructions are provided for both cases. In the subsequent section 1-$\gmb$-critical, 1-$\gmb'$-critical, and 2-$\gmb$-critical graphs are characterized, while 2-$\gmb'$-critical graphs are characterized among bipartite graphs. Finally, $\gsmb$- and $\gsmb'$-critical graphs are formally defined and characterized for graphs $G$ with $\gsmb(G),\gsmb'(G) \in \{1,2\}$  in Section~\ref{sec: SMBD-critical}.

\section{ Preliminaries}
\label{sec: prelim}

Let $G=(V(G), E(G))$ be a  graph. The order of  $G$ is denoted by $n(G)$. For a vertex $v\in V(G)$,  its open neighbourhood is denoted by $N(v)$ and its closed neighbourhood by $N[v]$.  The degree of $v$  is $\deg(v)=|N(v)|$.  The minimum and the maximum degree of $G$ are respectively denoted by $\delta(G)$ and $\Delta(G)$. An isolated vertex is a vertex of degree $0$, a leaf is a vertex of degree $1$.  A {\em support vertex} is a vertex adjacent to a leaf and a  {\em strong support vertex} is a vertex that is adjacent to at least two leaves. A set $S \subseteq V(G)$ is $x$-free for $x\in V(G)$ if $x\notin S$. As usual, $K_n$ and $P_n$ denote the complete graph and the path of order $n$, while $K_{m,n}$ is the complete bipartite graph of order $m+n$.

A {\em dominating set} of $G$ is a set $D \subseteq V(G)$ such that each vertex from $V(G)\setminus D$ has a neighbour in $D$. The {\em domination number} $\gamma(G)$ of $G$  is the minimum cardinality of a dominating set of $G$. If $X$ is a dominating set of $G$ with $|X| = \gamma(G)$, then $X$ is a {\em $\gamma$-set} of $G$. A vertex of degree $n(G) -1$ is a {\em dominating vertex}. If $G$ is a connected bipartite graph, then a vertex $x$ is a {\em bipartite dominating vertex} if $x$ is adjacent to all the vertices of the bipartition set of $G$ which does not contain $x$. An edge of a graph $G$ is a {\em dominating edge} if it is adjacent to all the other edges of $G$.   

The {\em outcome} $o(G)$ of the MBD game played on $G$ can be one of $\mathcal{D}$, $\mathcal{S}$, and $\mathcal{N}$, where $o(G)=\mathcal{D}$, if Dominator has a winning strategy no matter who starts the game; $o(G)=\mathcal{S}$, if Staller has a winning strategy no matter who starts the game; and $o(G)=\mathcal{N}$, if the first player has a winning strategy. See~\cite{duchene-2020} that the fourth possible option for the outcome never happens. We also add that in an optimal strategy of Dominator to achieve $\gmb(G)$ or $\gmb'(G)$, it is never an advantage for him to skip  a move. Moreover, if Staller skips a move it can never disadvantage Dominator~\cite[Lemma 2.3]{gledel-2019}. The same holds for the games in which Staller wins~\cite[Lemma 2.6]{bujtas-2024}.

In the rest of the preliminaries, we recall known results needed later. For $X\subseteq V(G)$, let $G|X$ denote the graph $G$ in which vertices from $X$ are considered as being already dominated. Then we have: 

\begin{theorem} [Continuation Principle~\cite{gledel-2019}]
\label{thm: Continuation Principle}
Let $G$ be a graph with $A, B \subseteq V(G)$. If $B\subseteq A$ then $\gmb(G|A) \leq \gmb(G|B)$ and  $\gmb'(G|A) \leq \gmb'(G|B)$.
\end{theorem}

\begin{proposition}{\rm \cite{bujtas-2024}} If $G$ is a graph, then the following properties hold. 
\begin{enumerate}
     \item If $o(G)=\mathcal{D}$ then $o(G+e)=\mathcal{D}$ for every $e\notin E(G)$.
	\item If $o(G)=\mathcal{S}$ then $o(G-e)=\mathcal{S}$ for every $e\in E(G)$.
	\item If $o(G)=\mathcal{N}$ then $o(G+e)\in \mathcal{\{N,D\}}$ for every $e\notin E(G)$.
	\item If $o(G)=\mathcal{N}$ then $o(G-e)\in \mathcal{\{N,S\}}$ for every $e\in E(G)$.
 \end{enumerate}
 \end{proposition}
\begin{lemma}{\rm \cite{dokyeesun-2024+}}\label{lem:doky}
If $G$ is a graph, then the following properties hold. 
 \begin{enumerate}
 \item[(i)] $\gmb(G)\leq \gmb(G-e)$ for every $e\in E(G)$.
\item[(ii)] $\gmb'(G)\leq \gmb'(G-e)$ for every $e \in E(G)$.
\item[(iii)] $\gsmb(G)\leq \gsmb(G+e)$ for every $ e \notin E(G)$.
\item[(iv)] $\gsmb'(G)\leq \gsmb'(G+e)$ for every $ e \notin E(G)$.  
\end{enumerate}
\end{lemma} 
 
\section{MBD game critical graphs}
\label{sec:main}

In this section, we introduce MBD game critical graphs. It is known from {\rm \cite{bujtas-2024}} that the outcome of the MBD game of a graph $G$ may change when an edge is removed or added. By Lemma~\ref{lem:doky}, the MBD number of a graph $G$ never decreases by removing an edge. This motivates us to define the MBD game critical graphs as follows.

\begin{definition}\label{definition 1} 
If $G$ is a graph and $\tau \in \{\gmb,\gmb'\}$, then $G$ is $\tau$-critical, if $\tau(G) < \tau(G-e)$, for any $e \in E(G)$. 
\end{definition}

Note that in view of Lemma~\ref{lem:doky}~(i) and (ii), if $G$ is $\tau$-critical, where $\tau \in \{\gmb,\gmb'\}$, then $\tau(G) < \infty$.  If $G$ is $\tau$-critical and $\tau(G)=k$, where $\tau \in \{\gmb,\gmb'\}$, then we say that $G$ is a {\em $k$-$\tau$-critical}. 

Let $\gmb(G)=k$. Iteratively removing edges $e$ with the property $\gmb(G-e)=k$ we arrive at a $k$-$\gmb$-critical spanning subgraph of $G$. Analogous conclusion holds for a $k$-$\gmb'$-critical graph. So such critical graphs clearly exists. In the rest of the section we prove that also connected critical graphs exist.

Consider graphs $G_k$, $k\geq 1$, constructed as follows. First, take the disjoint union of $k$ copies of $K_{2,2}$ with respective  bipartitions $\{x_i, x_i'\}$, $\{y_i, y_i'\}$, where $i\in[k]$. Then add vertex $w$ and make it adjacent to $x_i$ and $x_i'$ for $i\in [k]$. Finally, add a vertex $w'$ and the edge $ww'$. See Fig.~\ref{fig1:ilustrative-example}.

\begin{figure}[ht!]
\begin{center}
\begin{tikzpicture}[scale=0.5,style=thick,x=1cm,y=1cm]
\def\vr{3pt}
\begin{scope}[xshift=0cm, yshift=0cm] 
\coordinate(y_1) at (1,0);
\coordinate(y_1') at (3,0);
\coordinate(x_1) at (1,4);
\coordinate(x_1') at (3,4);
\draw (x_1) -- (y_1);
\draw (x_1) -- (y_1'); 
\draw (x_1') -- (y_1);
\draw (x_1') -- (y_1');
\draw (x_1) -- (10,8);
\draw (x_1') -- (10,8);
\draw(x_1)[fill=white] circle(\vr);
\draw(x_1')[fill=white] circle(\vr);
\draw(y_1)[fill=white] circle(\vr);
\draw(y_1')[fill=white] circle(\vr);
\node at (1,-.5) {$y_1$};
\node at (3,-.5) {$y_1'$};
\node at (0.3,4) {$x_1$};
\node at (2.3,4) {$x_1'$};
\end{scope}

\begin{scope}[xshift=4cm, yshift=0cm] 
\coordinate(y_2) at (1,0);
\coordinate(y_2') at (3,0);
\coordinate(x_2) at (1,4);
\coordinate(x_2') at (3,4);
\draw (x_2) -- (y_2);
\draw (x_2) -- (y_2'); 
\draw (x_2') -- (y_2);
\draw (x_2') -- (y_2');
\draw (x_2) -- (6,8);
\draw (x_2') -- (6,8);
\draw(x_2)[fill=white] circle(\vr);
\draw(x_2')[fill=white] circle(\vr);
\draw(y_2)[fill=white] circle(\vr);
\draw(y_2')[fill=white] circle(\vr);
\node at (1,-0.5) {$y_2$};
\node at (3,-0.5) {$y_2'$};
\node at (.4,4) {$x_2$};
\node at (2.35,4) {$x_2'$};
\end{scope}

\begin{scope}[xshift=8cm, yshift=0cm] 

\coordinate(w) at (2,8);

\coordinate(w') at (2,12);
\draw (w) -- (w');

\draw(w)[fill=white] circle(\vr);
\draw(w')[fill=white] circle(\vr);

\node at (2.5,8.15) {$w$};
\node at (2,12.5) {$w'$};

\node at (1.5,0) {$\cdots$};
\node at (1.5,2) {$\cdots$};
\node at (1.5,4) {$\cdots$};
\end{scope}

\begin{scope}[xshift=12cm, yshift=0cm] 
\coordinate(y_{k-1}) at (1,0);
\coordinate(y_{k-1}') at (3,0);
\coordinate(x_{k-1}) at (1,4);
\coordinate(x_{k-1}') at (3,4);
\draw (x_{k-1}) -- (y_{k-1});
\draw (x_{k-1}) -- (y_{k-1}'); 
\draw (x_{k-1}') -- (y_{k-1});
\draw (x_{k-1}') -- (y_{k-1}');
\draw (x_{k-1}') -- (-2,8);
\draw (x_{k-1}) -- (-2,8);
\draw(x_{k-1})[fill=white] circle(\vr);
\draw(x_{k-1}')[fill=white] circle(\vr);
\draw(y_{k-1})[fill=white] circle(\vr);
\draw(y_{k-1}')[fill=white] circle(\vr);
\node at (1,-0.5) {$y_{k-1}$};
\node at (3,-0.5) {$y_{k-1}'$};
\node at (.1,4) {$x_{k-1}$};
\node at (2.1,4) {$x_{k-1}'$};
\end{scope}

\begin{scope}[xshift=16cm, yshift=0cm] 
\coordinate(y_k) at (1,0);
\coordinate(y_k') at (3,0);
\coordinate(x_k) at (1,4);
\coordinate(x_k') at (3,4);
\draw (x_k) -- (y_k);
\draw (x_k) -- (y_k'); 
\draw (x_k') -- (y_k);
\draw (x_k') -- (y_k');
\draw (x_k') -- (-6,8);
\draw (x_k) -- (-6,8);
\draw(x_k)[fill=white] circle(\vr);
\draw(x_k')[fill=white] circle(\vr);
\draw(y_k)[fill=white] circle(\vr);
\draw(y_k')[fill=white] circle(\vr);
\node at (1,-0.5) {$y_k$};
\node at (3,-0.5) {$y_k'$};
\node at (1.8,4.) {$x_k$};
\node at (3.8,4) {$x_k'$};
\end{scope}

\end{tikzpicture}
\caption{The graph $G_k$ which is a connected $(k+1)$-$\gmb$-critical graph}
\label{fig1:ilustrative-example}
\end{center}
\end{figure}

\begin{proposition}
    If $k\geq 1$, then $G_k$ is a $(k+1)$-$\gmb$-critical graph.
\end{proposition}
\proof
We first note that $\gamma(G_k)=k+1$, so that $\gmb(G_k)\geq k+1$. Assume now that Dominator starts the game on $G_k$ by selecting the vertex $w$. Then in the rest of the game, he can select one of the vertices $x_i$ and $x_i'$ for each $i\in[k]$. It follows that $\gmb(G_k)\leq k+1$.

To show that $\gmb(G_k-e)>k+1$ for any edge $e\in E(G_k)$, by the symmetry of $G_k$ it suffices to consider three typical edges. Let first $e=ww'$. Then Dominator must start the game by playing $w'$. Since the domination number of the large component of $G_k-ww'$ is $k+1$, it follows that 
$\gmb(G_k-ww')\geq k+2$. Consider next the edge $e=wx_1$. Clearly, $w$ is an optimal first move of Dominator. Now Staller replies with the move $x_1$. Since $x_1$, $y_1$, and $y_1'$ are not yet dominated, Dominator will need two moves to dominate them. This in turn implies that $\gmb(G_k-wx_1)\geq k+2$. Consider finally the edge $e=x_1y_1$. In this case, we again see that $w$ is an optimal first move of Dominator and if now Staller replies by playing $x_1'$ we can see as in the previous case that $\gmb(G_k-x_1y_1)\geq k+2$.
\qed

To show that there exist connected $k$-$\gmb'$-critical graphs, consider graphs $H_k$, $k\geq 1$, obtained as follows. First, take the disjoint union of $k+1$ copies of $K_{2,3}$ whose respective   bipartitions are $\{x_i, x_i'\}$, $\{y_i, y_i', y_i''\}$, where $i\in[k+1]$. Then add all possible edges between $x_{k+1}$ and $x_1, x_1',\ldots , x_k,x_k'$, and between $x_{k+1}'$ and $x_1, x_1',\ldots , x_k,x_k'$. See Fig.~\ref{fig2:ilustrative-example}.

\begin{figure}[ht!]
\begin{center}
\begin{tikzpicture}[scale=0.5,style=thick,x=1cm,y=1cm]
\def\vr{3pt}
\begin{scope}[xshift=0cm, yshift=0cm] 
\coordinate(y_1) at (.5,0);
\coordinate(y_1') at (2,0);
\coordinate(y_1'') at (3.5,0);
\coordinate(x_1) at (.5,4);
\coordinate(x_1') at (3.5,4);
\draw (x_1) -- (y_1);
\draw (x_1) -- (y_1'); 
\draw (x_1) -- (y_1'');
\draw (x_1') -- (y_1);
\draw (x_1') -- (y_1');
\draw (x_1') -- (y_1'');
\draw (x_1) -- (8.5,8);
\draw (x_1') -- (8.5,8);
\draw (x_1) -- (11.5,8);
\draw (x_1') -- (11.5,8);

\draw(x_1)[fill=white] circle(\vr);
\draw(x_1')[fill=white] circle(\vr);
\draw(y_1)[fill=white] circle(\vr);
\draw(y_1')[fill=white] circle(\vr);
\draw(y_1'')[fill=white] circle(\vr);
\node at (0,-.5) {$y_1$};
\node at (2,-.5) {$y_1'$};
\node at (3.5,-.5) {$y_1''$};
\node at (-.5,4) {$x_1$};
\node at (2.5,4) {$x_1'$};

\end{scope}

\begin{scope}[xshift=4cm, yshift=0cm] 
\coordinate(y_2) at (.5,0);
\coordinate(y_2') at (2,0);
\coordinate(y_2'') at (3.5,0);
\coordinate(x_2) at (.5,4);
\coordinate(x_2') at (3.5,4);
\draw (x_2) -- (y_2);
\draw (x_2) -- (y_2'); 
\draw (x_2) -- (y_2'');
\draw (x_2') -- (y_2);
\draw (x_2') -- (y_2');
\draw (x_2') -- (y_2'');
\draw (x_2) -- (4.5,8);
\draw (x_2') -- (4.5,8);
\draw (x_2) -- (7.5,8);
\draw (x_2') -- (7.5,8);
\draw(x_2)[fill=white] circle(\vr);
\draw(x_2')[fill=white] circle(\vr);
\draw(y_2)[fill=white] circle(\vr);
\draw(y_2')[fill=white] circle(\vr);
\draw(y_2'')[fill=white] circle(\vr);
\node at (.5,-.5) {$y_2$};
\node at (2,-.5) {$y_2'$};
\node at (4,-.5) {$y_2''$};
\node at (1.5,4) {$x_2$};
\node at (4.5,4) {$x_2'$};
\end{scope}

\begin{scope}[xshift=8cm, yshift=0cm] 
\coordinate(y_{k+1}) at (.5,12);
\coordinate(y_{k+1}') at (2,12);
\coordinate(y_{k+1}'') at (3.5,12);
\coordinate(x_{k+1}) at (.5,8);
\coordinate(x_{k+1}') at (3.5,8);
\draw (x_{k+1}) -- (y_{k+1});
\draw (x_{k+1}) -- (y_{k+1}'); 
\draw (x_{k+1}) -- (y_{k+1}'');
\draw (x_{k+1}') -- (y_{k+1});
\draw (x_{k+1}') -- (y_{k+1}');
\draw (x_{k+1}') -- (y_{k+1}'');

\draw(x_{k+1})[fill=white] circle(\vr);
\draw(x_{k+1}')[fill=white] circle(\vr);
\draw(y_{k+1})[fill=white] circle(\vr);
\draw(y_{k+1}')[fill=white] circle(\vr);
\draw(y_{k+1}'')[fill=white] circle(\vr);
\node at (-.2,12.7) {$y_{k+1}$};
\node at (1.7,12.7) {$y_{k+1}'$};
\node at (3.6,12.7) {$y_{k+1}''$};
\node at (-.5,8.2) {$x_{k+1}$};
\node at (4.9,8.2) {$x_{k+1}'$};
\node at (1.8,0) {$\cdots$};
\node at (1.8,2) {$\cdots$};
\node at (1.8,4) {$\cdots$};
\end{scope}

\begin{scope}[xshift=12cm, yshift=0cm] 
\coordinate(y_{k-1}) at (.5,0);
\coordinate(y_{k-1}') at (2,0);
\coordinate(y_{k-1}'') at (3.5,0);
\coordinate(x_{k-1}) at (.5,4);
\coordinate(x_{k-1}') at (3.5,4);
\draw (x_{k-1}) -- (y_{k-1});
\draw (x_{k-1}) -- (y_{k-1}'); 
\draw (x_{k-1}) -- (y_{k-1}'');
\draw (x_{k-1}') -- (y_{k-1});
\draw (x_{k-1}') -- (y_{k-1}');
\draw (x_{k-1}') -- (y_{k-1}'');
\draw (x_{k-1}) -- (-3.5,8);
\draw (x_{k-1}') -- (-3.5,8);
\draw (x_{k-1}) -- (-.5,8);
\draw (x_{k-1}') -- (-.5,8);
\draw(x_{k-1})[fill=white] circle(\vr);
\draw(x_{k-1}')[fill=white] circle(\vr);
\draw(y_{k-1})[fill=white] circle(\vr);
\draw(y_{k-1}')[fill=white] circle(\vr);
\draw(y_{k-1}'')[fill=white] circle(\vr);
\node at (0,-.5) {$y_{k-1}$};
\node at (2,-.5) {$y_{k-1}'$};
\node at (3.7,-.55) {$y_{k-1}''$};
\node at (-.5,4) {$x_{k-1}$};
\node at (2.5,4) {$x_{k-1}'$};

\end{scope}

\begin{scope}[xshift=16cm, yshift=0cm] 
\coordinate(y_k) at (.5,0);
\coordinate(y_k') at (2,0);
\coordinate(y_k'') at (3.5,0);
\coordinate(x_k) at (.5,4);
\coordinate(x_k') at (3.5,4);
\draw (x_k) -- (y_k);
\draw (x_k) -- (y_k'); 
\draw (x_k) -- (y_k'');
\draw (x_k') -- (y_k);
\draw (x_k') -- (y_k');
\draw (x_k') -- (y_k'');
\draw (x_k) -- (-7.5,8);
\draw (x_k') -- (-7.5,8);
\draw (x_k) -- (-4.5,8);
\draw (x_k') -- (-4.5,8);
\draw(x_k)[fill=white] circle(\vr);
\draw(x_k')[fill=white] circle(\vr);
\draw(y_k)[fill=white] circle(\vr);
\draw(y_k')[fill=white] circle(\vr);
\draw(y_k'')[fill=white] circle(\vr);
\node at (.8,-.5) {$y_k$};
\node at (2.2,-.5) {$y_k'$};
\node at (4,-.5) {$y_k''$};
\node at (1.5,4) {$x_k$};
\node at (4.5,4) {$x_k'$};
\end{scope}

\end{tikzpicture}
\caption{The graph $H_k$ which is a connected $(k+1)$-$\gmb'$-critical graph}
\label{fig2:ilustrative-example}
\end{center}
\end{figure}

\begin{proposition}
    If $k\geq 1$, then $H_k$ is a $(k+1)$-$\gmb'$-critical graph.
\end{proposition}
\proof
Since $\gamma(H_k)=k+1$ and because during the S-game Dominator is able to select one vertex from each of the sets $\{x_i,x_i'\}$, $i\in[k+1]$, we infer that $\gmb'(H_k)=k+1$.

To show that $H_k$ is $(k+1)$-$\gmb'$-critical, by the symmetry of $H_k$ it suffices to consider three typical edges of $H_k$.
Assume first that $e=x_{k+1}y_{k+1}$. Then in the S-game played on $H_k-x_{k+1}y_{k+1}$, Staller's  strategy is that she first selects the vertex $x_{k+1}'$. Then Dominator must reply by choosing the vertex $y_{k+1}$ for otherwise Staller wins in her next move by selecting $y_{k+1}$. Then  Staller selects $x_{k+1}$ as her second move. Because $y'_{k+1}$ and $y''_{k+1}$ are adjacent only to both $x_{k+1}$ and $x'_{k+1}$, Staller will be able to win the game in her third move by either selecting $y'_{k+1}$ or $y''_{k+1}$. As a consequence, $\gmb'(H_k-x_{k+1}y_{k+1})=\infty$. 
The next typical edge to be considered is $e=x_1x_{k+1}$. Then Staller in her strategy first selects the vertex $x_{k+1}'$. This move forces Dominator to play $x_{k+1}$. As her second move, Staller then plays $x_1$  which in turn forces Dominator to play $x_1'$. But now since $x_1$ is not yet dominated, we can conclude that in the rest of the game, Dominator must play at least $k$ more moves. Hence also in this case, at least $k+2$  vertices will be selected by him.
The last typical edge to be considered is the edge $x_1y_1$. In this case, Staller first selects $x_1'$ which forces Dominator to play $y_1$ (otherwise Staller will win in her next move). Now Staller plays $x_1$ and then she wins in her next move. Hence $\gmb'(H_k-x_1y_1)=\infty$.  
\qed

\section{Critical graphs with small MBD numbers}
\label{sec:MBD-critical}

In this section, we describe critical graphs for the cases in which Dominator wins the game in one or two moves.

\begin{proposition}\label{thm mbd1}
A connected graph $G$ is  $1$-$\gmb$-critical  if and only if $G=K_{1,n}, n \geq 1$.
\end{proposition}
\proof 
First, assume that $G$ is connected and  $1$-$\gmb$-critical. Then $G$ contains a dominating vertex, say $u$. If possible, suppose that there exists an edge $e$ of $G$ that is not incident with $u$. Then $\gmb(G-e) = 1$, which leads to a contradiction. Hence every edge of $G$ is incident with $u$. Thus $G$ is isomorphic to $K_{1,n}, n \geq 1$.

Conversely, let $G$ be a star $K_{1,n}$ for $n\geq 1$.  Clearly, $G$ is connected and $\gmb(G)=1$.  Deletion of any edge of $G$ results in a disconnected graph and hence Dominator cannot win this game in one move. Therefore $G$ is $1$-$\gmb$-critical.
\qed 

Let $K'_{2,n}$, $n\ge 1$, be the complete split graph (cf.~\cite{ge-2024}) consisting of a clique of order $2$ and an independent set of order $n$, where every vertex in the independent set is adjacent to both vertices of the clique. (We use this notation because $K'_{2,n}$ can be obtained from $K_{2,n}$ by adding a single edge.) 

\begin{theorem}\label{thm mbd2}
A connected graph $G$ is $1$-$\gmb'$-critical  if and only if $G = K'_{2,n}$, $n\ge 1$.    
\end{theorem}

\proof
Assume that $G$ is a $1$-$\gmb'$-critical graph.  Therefore, Dominator can win this game by selecting one vertex in his first move as a second player. Thus $G$ has at least two dominating vertices say $u$ and $v$. Now we show that $e=uv$ is a dominating edge. Suppose on the contrary that there exists an edge $f=xy$ such that $x,y \notin \{u,v\}$. Then $\gmb'(G-f)=1$, which is a contradiction. Therefore, every edge of $G$ is incident with either $u$ or $v$  or both. Hence the edge $e=uv$ is a dominating edge of  $G$. Since $u$ and $v$ are dominating vertices, every vertex other than $u$ and $v$ has degree two in $G$. It follows that $G = K'_{2,n}$ for some $n\ge 1$.     

Conversely, assume that $G = K'_{2,n}$ for some $n\ge 1$. Then $G$ has two dominating vertices and hence $\gmb'(G)=1$. Moreover, $G-e$ has at most one dominating vertex for every $e\in E(G)$. Thus in the S-game played on $G-e$, Dominator needs at least two moves because Staller can choose the dominating vertex of $G-e$ in her first move if there is such a vertex. Therefore $G$ is  $1$-$\gmb'$-critical.
\qed
 
\begin{theorem}\label{thm mbd3}
A connected graph $G$ is  $2$-$\gmb$-critical   if and only if $G$ is obtained from a star $K_{1,n}$, $n\geq 1$, with center $u$ and a $K_{2,m}$, $m\geq 2$, whose  bipartition is $\{x_1, x_2\}$, $\{y_1, y_2,\ldots, y_m\}$,  by adding  the edges $ux_1$ and $ux_2$ (see Fig.~\ref{fig3:ilustrative-example}).
\end{theorem}

\begin{figure}[ht!]
\begin{center}
\begin{tikzpicture}[scale=.5,style=thick,x=1cm,y=1cm]
\def\vr{3pt}
\begin{scope}[xshift=0cm, yshift=0cm] 
\coordinate(u) at (2.5,6);
\coordinate(u_1) at (-3.5,9);
\coordinate(u_2) at (-2,9);
\coordinate(u_3) at (-.5,9);
\coordinate(u_4) at (1,9);
\coordinate(u_5) at (2.5,9);
\coordinate(u_6) at (4,9);
\coordinate(u_7) at (5.5,9);
\coordinate(u_8) at (7,9);
\coordinate(u_9) at (8.5,9);

\coordinate(x_1) at (1,4);
\coordinate(x_2) at (4,4);

\coordinate(y_1) at (-4.5,0);
\coordinate(y_2) at (-3,0);
\coordinate(y_3) at (-1.5,0);
\coordinate(y_4) at (0,0);
\coordinate(y_5) at (1.5,0);
\coordinate(y_6) at (3,0);
\coordinate(y_7) at (4.5,0);
\coordinate(y_8) at (6,0);
\coordinate(y_9) at (7.5,0);
\coordinate(y_10) at (9,0);
\coordinate(y_11) at (10.5,0);

\draw (x_1) -- (u);
\draw (x_2) -- (u);

\draw (u) -- (u_1);
\draw (u) -- (u_2);
\draw (u) -- (u_3);
\draw (u) -- (u_7);
\draw (u) -- (u_8);
\draw (u) -- (u_9);

\draw (x_1) -- (y_1);
\draw (x_1) -- (y_2);
\draw (x_1) -- (y_3);
\draw (x_1) -- (y_8);
\draw (x_1) -- (y_9);
\draw (x_1) -- (y_10);
\draw (x_1) -- (y_11);

\draw (x_2) -- (y_1);
\draw (x_2) -- (y_2);
\draw (x_2) -- (y_3);
\draw (x_2) -- (y_8);
\draw (x_2) -- (y_9);
\draw (x_2) -- (y_10);
\draw (x_2) -- (y_11);

\draw(u)[fill=white] circle(\vr);
\draw(x_1)[fill=white] circle(\vr);
\draw(x_2)[fill=white] circle(\vr);

\foreach \i in {1,...,3}
{ \draw(u_\i)[fill=white] circle(\vr);}
\foreach \i in {7,...,9}
{ \draw(u_\i)[fill=white] circle(\vr);}

\foreach \i in {1,...,3}
{ \draw(y_\i)[fill=white] circle(\vr);}
\foreach \i in {8,...,11}
{ \draw(y_\i)[fill=white] circle(\vr);}

\node at (2.5,0) {$\cdots$};
\node at (2.5,9) {$\cdots$};

\node at (3.1,5.9) {$u$};
\node at (-3.5,9.5) {$u_1$};
\node at (-2,9.5) {$u_2$};
\node at (-0.5,9.5) {$u_3$};
\node at (5.45,9.5) {$u_{n-2}$};
\node at (7.2,9.5) {$u_{n-1}$};
\node at (8.6,9.5) {$u_n$};

\node at (.4,4.1) {$x_1$};
\node at (4.6,4.1) {$x_2$};

\node at (-4.5,-.5) {$y_1$};
\node at (-3,-.5) {$y_2$};
\node at (-1.5,-.5) {$y_3$};
\node at (7.4,-.5) {$y_{m-2}$};
\node at (9.7,-.5) {$y_{m-1}$};
\node at (11.2,-.5) {$y_m$};
\end{scope}
\end{tikzpicture}

\caption{Connected $2$-$\gmb$-critical graphs}
\label{fig3:ilustrative-example}
\end{center}
\end{figure}
\proof
Assume that $G$ is connected and 2-$\gmb$-critical. Since $\gmb(G)=2$, Dominator cannot finish the game in one move and hence $G$ has no dominating vertex. 

Let $u$ be an optimal first move of Dominator in the D-game played on $G$. Suppose that there exists an edge $e$ with both end-vertices in $N(u)$. Then, Dominator can win the D-game on $G-e$ in two moves using the same strategy as that of the D-game on $G$  which is a contradiction with the fact that $G$ is a 2-$\gmb$-critical graph. Thus $N(u)$ is an independent set. 

Further, we prove that $|V(G) \setminus N[u]| \geq 2$. Since $u$ is not a dominating vertex, $|V(G) \setminus N[u]| \geq 1$. If possible suppose that there exists exactly one vertex $x\in V(G)\setminus N[u]$. Since $G$ is connected, $x$ is adjacent to at least one vertex, say $w$, in $N(u)$. Let $e=uw$ be an edge of $G$ and consider a D-game on $G-e$. Dominator selects $u$ as his first move in $G-e$. The only undominated vertices of $G$ are $x$  and $w$. Therefore, Dominator can finish this game in his next move by selecting either $x$ or $w$ depending on the Staller's move and this contradicts that $G$ is 2-$\gmb$-critical.

Since Dominator has a winning strategy with two moves, there exists $a,b \in V(G)\setminus \{u\}$ that are both adjacent to all vertices from $V(G)\setminus N(u)$. If both $a,b$ are from $V(G)\setminus N(u)$, then $\gmb(G-e)=2$ holds for any edge $e$ between $N(u)$ and  $V(G) \setminus N[u]$ (at least one such edge must exist, since $G$ is connected), which contradicts the fact that $G$ is 2-$\gmb$-critical. Thus assume that $a \in N(u)$. Since $\gmb(G-ua) > 2$,  Dominator cannot finish the D-game played on $G-ua$ by first selecting $u$ and then in his second move one vertex from $\{a,b\}$ (that was not selected by Staller in her first move). This is possible only if Staller in her first move selects $a$ and $ab \notin E(G)$. Since $a$ dominates all vertices of $V(G) \setminus N[u]$, this implies that $b \in N(u)$. Hence all vertices of $G$ that dominate the whole $G-N(u)$ are from $N(u)$. Hence in the D-game played on $G$ Dominator's optimal second move will be a vertex from $N(u)$.

Let $a$ and $b$ be two vertices of $N(u)$ that dominate all vertices in $V(G)\setminus N[u]$. Suppose that there exists a vertex $x \in V(G) \setminus N[u]$ that is adjacent to $y \in V(G)\setminus \{a,b\}$. Then $\gmb(G-xy)=2$, a contradiction. Hence vertices in $V(G) \setminus N[u]$ have degree 2 in $G$ and are adjacent to $a$ and $b$.

Finally, we show that $|N(u)| \geq 3$. If possible suppose that $a$ and $b$ are the only neighbours of $u$ in $G$. Let $e=ua$. Then $\gmb(G-e)=2$. Indeed, Dominator can select $b$ as his first optimal move in a D-game played on $G-e$. The only vertex undominated after this move is $a$. Therefore Dominator can finish the game by selecting a vertex in $V(G) \setminus N[u]$ as his next move. Hence we again get a contradiction. Thus $|N[u]| \geq 3$.

By the above properties of the graph $G$ we can deduce that vertices in $\{u\} \cup (N(u)\setminus \{a,b\})$ induce $K_{1,n}$ for some $n \geq 1$, vertices in the set $\{a,b\} \cup (V(G) \setminus N[u])$ induce $K_{2,m}$ for $m\geq 2$ and $au,bu$ are the only edges between $K_{1,n}$ and $K_{2,m}$ in $G$.

Conversely, let $G$ be a graph obtained from $K_{1,n}$, $n\geq 1$, with center $u$, and from $K_{2,m}$, $m\geq 2$, with bipartition $\{x_1, x_2\}$, $\{y_1, y_2,\ldots, y_m\}$, by adding  the edges $ux_1$ and $ux_2$. Clearly, $G$ is connected and Dominator can finish a D-game on $G$ by selecting $u$ as his first move and then selecting either $x_1$ or $x_2$ with respect to the Staller's first move. So $\gmb(G)=2$. 

Let $e$ be a pendant edge incident to $u$. Clearly, the graph $G-e$ has an isolated vertex. Dominator selects this isolated vertex as his first optimal move in a D-game on $G-e$. And the remaining part of $G-e$  has no dominating vertices. Therefore Dominator needs at least three moves to finish a D-game on $G-e$.

Now let $e=ux_1$. Consider a D-game on $G-e$. Since $u$ is a support vertex Domiantor first selects $u$. Then Staller selects $x_1$. If Dominator selects $x_2$ as his next move then $x_1$ remains undominated. And if Dominator selects a vertex in $V(G) \setminus N[u]$, then there is an undominated vertex in $V(G) \setminus N[u]$. Therefore Dominator needs at least three moves to finish the game in $G-e$. A similar argument also holds for $e=ux_2$.

Now let $e$ be an edge whose one end vertex is $x_1$ and the other end vertex lies in $V(G) \setminus N[u]$. Consider a D-game on $G-e$. Clearly, $u$ and $x_2$ are support vertices. Definitely, Staller can select one of these support vertices in her turn. Therefore Dominator must select the pendant vertex adjacent to the support vertex selected by Staller. Clearly, this restriction does not allow Dominator to finish the game on $G-e$ in two moves. Hence $\gmb(G-e) > 2$ in this case and this is the same when $e$ is an edge between $x_2$ and $V(G) \setminus N[u]$. 

From all these cases we conclude that $G$ is 2-$\gmb$-critical.
\qed

For $2$-$\gmb'$-critical graphs we do not have a complete characterization but can give the following necessary conditions. 

\begin{proposition} \label{Prop:necessary-conditions}
If $G$ is connected $2$-$\gmb'$-critical graph, then the following properties hold. 
\begin{enumerate}
    \item[(i)] $n(G)\ge 5$.
    \item[(ii)] $\delta(G) \geq 2$.
    \item[(iii)] $\Delta(G) \leq n(G)-2$.
\end{enumerate}
Moreover, all the bounds are sharp. 
\end{proposition}

\proof
(i) Let $G$ be a connected $2$-$\gmb'$-critical graph. Since $\gmb'(G)=2$, $G$ has at least four vertices (note that Staller has two moves and Dominator has two moves). Suppose that $G$ is 2-$\gmb'$-critical connected graph with $n(G)=4$. Then $G$ contains $2K_2$ as strict spanning subgraph and thus $G$ cannot be 2-$\gmb'$-critical. 

\medskip
(ii) Suppose on the contrary that $\delta(G)=1$. Let $u \in V(G)$ such that $\deg (u)=1$  and let $v$ be the only neighbor of $u$. Since $G$ is connected and has at least $5$ vertices by (i), we have $\deg (v)\geq 2$.

Now consider an S-game on $G$. The vertex $v$ is an optimal first move of Staller, cf.~\cite{duchene-2020}. Therefore Dominator must play $u$ as his first reply, otherwise, Staller will win this game by selecting $u$.  Since $\gmb'(G)=2$, Dominator can dominate all the vertices of $G$ other than $u$ and $v$ in his next move.  Hence there exist $a,b \in V(G)\setminus \{u,v\}$ that dominate all the vertices of $V(G) \setminus \{u,v\}$ or, equivalently, the subgraph of $G$ induced by $V(G) \setminus \{u,v\}$ has two dominating vertices. Let $e$ be an arbitrary edge incident with $v$ but not with $u$. Then $\gmb'(G-e)=2$, a contradiction with $G$ being 2-$\gmb'$-critical. Thus we conclude that $\delta(G)\geq 2$.

\medskip
(iii) Let $G$ be connected and $2$-$\gmb'$-critical.  Since $\gmb'(G) \neq 1$, $G$ has at most one dominating vertex. For the purpose of contradiction assume that $G$ contains a dominating vertex $u$, i.e.\ $\deg(u)=n(G)-1$. 
It is known from~\cite[Proposition~4.2]{divakaran-2024} that $\gmb'(G)=2$ implies that there exists a vertex $v\neq u$  such that  $\{v, v_1\}$ and  $\{v,v_2\}$ are two $u$-free $\gamma$-sets. Since $\deg(v)< n-1$,  there is a vertex $w$ in $G$  which is not adjacent to $v$. We consider two cases for the remaining part.

\medskip\noindent 
\textbf{Case 1}: $w\notin \{v_1,v_2\}$.\\
Since $\{v, v_1\}$ and $\{v,v_2\}$  are two $u$-free $\gamma$-sets of $G$, the vertex $w$ must be dominated by both $v_1$ and $v_2$. Let $e=uw$. Consider an S-game on $G-e$. If \textit{}Staller first selects $u$, then Dominator can select all the vertices from either
$\{v, v_1\}$ or $\{v,v_2\}$ and win the game in two moves.
    
If Staller first selects a vertex other than $u$, then Dominator must select $u$.  In this case, the only undominated vertex is $w$.  Since $w$ is dominated by all the vertices from the set $A=\{w, v_1, v_2\}$, Dominator can select any vertex of $A$ (depending on Staller's move) and win the game in two moves.
This contradicts that $G$ is $2$-$\gmb'$-critical.

\medskip\noindent 
\textbf{Case 2}: $w\in \{v_1,v_2\}$.\\
In this case, $v_1v_2$ must be an edge of $G$. Let $e=uv_1$. Consider an S-game on $G-e$. If Staller selects a vertex from the set $\{v_1, v_2\}$, then Dominator selects the other one. Now Dominator selects one vertex from $\{u,v\}$ (depending on the Staller's move)  and wins the game in his second move. This leads to a contradiction to the assumption that $G$ is $2$-$\gmb'$-critical. Now, if Staller selects $u$, then Dominator can select all the vertices from either  $\{v, v_1\}$ or  $\{v,v_2\}$  in his first two moves and finishes the game in two moves, a contradiction.  

\medskip
To prove that all three bounds are sharp, consider the complete bipartite graphs $K_{2,n}$, $n\geq 3$, which are connected and $2$-$\gmb'$-critical. 
\qed

To conclude the section we characterize connected bipartite graphs that are $2$-$\gmb'$-critical.The result will follow from the following two lemmas.

\begin{lemma} \label{lem:mbd 7}
Let $G$ be a bipartite graph with bipartition $V_1, V_2$, where $|V_i|\geq 3$, $i\in [2]$. If $G$ has exactly two bipartite dominating vertices in each $V_i$, $i\in [2]$, and every vertex has degree two except bipartite dominating vertices, then $G$ is  $2$-$\gmb'$-critical.
\end{lemma} 

\proof
Let $G$ be a graph satisfying the above properties. Let $v_{i,1}$ and $v_{i,2}$ be the  bipartite dominating vertices of $V_{i}$, $i\in[2]$. Clearly, 
$\{v_{1,1}, v_{2,1}\}$, and $\{v_{1,1}, v_{2,2}\}$  are two $v_{1,2}$-free  $\gamma$-sets of cardinality 2. Also $\{v_{1,2}, v_{2,1}\}$, and $\{v_{1,2}, v_{2,2}\}$  are two $v_{1,1}$-free  $\gamma$ sets of cardinality 2. Therefore Dominator can select two vertices from one set irrespective of Staller's move and win the game in two moves. Thus $\gmb'(G)=2$.  Now we show that $G$ is critical.
Any edge of $G$ either has two bipartite dominating vertices as endpoints or has exactly one bipartite dominating vertex as an endpoint. 
Let first $e=v_{1,1}v_{2,1}$. Consider an S-game on $G-e$. Clearly, the only bipartite dominating vertices in $G-e$ are $v_{1,2}$ and $v_{2,2}$.   In her strategy Staller first selects one of the bipartite dominating vertices, say $v_{1,2}$. If the first optimal move of Dominator is in $V_2$, then by the Theorem~\ref{thm: Continuation Principle}, $v_{2,2}$ is an optimal first move of Dominator.   Now all the vertices in $V_2$ except $v_{2,2}$ are undominated. Thus Staller selects $v_{1,1}$ as her next move. Any vertex in $V_1\setminus \{v_{1,1},v_{1,2}\}$ dominates only $v_{2,1}$ and $v_{2,2}$ in $V_2$. But $V_2$ has at least three vertices and Dominator cannot finish the game in two moves.  Now assume that an optimal first move of Dominator is in $V_1$. Let the first optimal move of Dominator be $v_{1,1}$ after the same first move of Staller. Now Staller selects $v_{2,1}$. So the undominated vertices are $v_{2,1}$ and all vertices in $V_1$ except $v_{1,1}$. Since $G$ is bipartite, Dominator needs at least two more moves to finish the game in this case.  Finally, assume that Dominator selects an unplayed vertex other than $v_{1,1}$ in $V_1$ after the same first move of Staller. This vertex only Dominates itself and both the vertices $v_{2,1}$  and $v_{2,2}$. It is clear that there are still undominated vertices in both $V_1$ and $V_2$. Therefore Dominator needs at least three moves to finish the game.  Hence we can conclude that $\gmb'(G-e)>2$.

If $e$ is an edge between any two bipartite dominating vertices, then $\gmb'(G-e)>2$ is proved by similar arguments as above.

Finally, let $e=ab$ be an edge, where exactly one of its end vertices, say $a$, is a bipartite dominating vertex.   Clearly, $b$ is the leaf of $G-e$. Therefore Staller first selects the support vertex of $G-e$ and then Dominator must select the leaf $b$ of  $G-e$  as his first move. We can see that there are still undominated vertices in both $V_1$ and $V_2$. So Dominator needs at least two more moves to finish the game. Thus we can conclude that $\gmb'(G-e)>2$. 

Therefore $G$ is a connected $2$-$\gmb'$-critical graph.
\qed

\begin{lemma}\label{lem:Bipartite2}
Let $G$ be a connected bipartite graph with bipartition $V_1,V_2$. If $G$ is $2$-$\gmb'$-critical, then $G$ is either $K_{2,m}$ for $m \geq 3$, or $|V_i| \geq 3$ for $i \in [2]$ and $G$ has exactly two bipartite dominating vertices in each $V_i$, $i \in [2]$ and all other vertices of $G$ are of degree $2$.
\end{lemma}

\proof
Let $V_1=\{v_{1,1},\ldots , v_{1,m}\}$ and $V_2=\{v_{2,1},\ldots , v_{2,n}\}$. First assume that $G$ is $2$-$\gmb'$-critical. If $|V_i|=1$ for some $i \in [2]$, then Proposition~\ref{Prop:necessary-conditions} implies that $|V_j| \geq 4$ for $j \in [2]\setminus \{i\}$ and thus $G=K_{1,\ell}$, $\ell \geq 4$. Hence $\gmb'(G) = \infty$, a contradiction with $G$ being $2$-$\gmb'$-critical. Therefore, since $n(G) \geq 5$ by Proposition~\ref{Prop:necessary-conditions}, we may without loss of generality assume that $|V_1|\geq 2$ and  $|V_2|\geq 3$. 

First assume that $|V_i|=2$ for some $i \in [2]$ and consequently $|V_j| \geq 3$ for $j \in [2]\setminus \{i\}$. Since $\delta(G) \geq 2$ by Proposition~\ref{Prop:necessary-conditions}, any vertex $x \in V_j$ is adjacent to both vertices of $V_i$. Hence $G$ is isomorphic to $K_{2,m}$, where $m =|V_j|\geq 3$. In the rest of the proof we may thus assume that $|V_i| \geq 3$ for $i \in [2]$.

Suppose that $V_i$ has at most one bipartite dominating vertex, say $v_{i,1}$ for some $i \in [2]$. Then Staller selects $v_{i,1}$ in her first move and she is able to selects her second vertex so that Dominator cannot finish the game within two moves, a contradiction. Therefore, each $V_i$ has at least two bipartite dominating vertices.

Suppose that $V_i$ for some $i\in [2]$ contains three bipartite dominating vertices, say $v_{i,1}, v_{i,2}$, and $v_{i,3}$. Let $v_{j,1}$ and $v_{j,2}$ be two bipartite dominating vertices of $V_j$. Consider the edge $e=v_{i,3}v_{j,3}$. Then $v_{i,1}, v_{j,1}, v_{i,2}, v_{j,2}$ are bipartite dominating vertices in $G-e$. Therefore Dominator has a strategy to  win the game on $G-e$ in two moves by selecting one of the bipartite dominating vertices in each $V_i$, $i\in [2]$, no matter how Staller plays. Thus $\gmb'(G-e)=2$ and hence $G$ is not $2$-$\gmb'$-critical, a contradiction. Therefore each $V_i$, $i\in [2]$, has exactly two bipartite dominating vertices, say $v_{1,1}, v_{2,1}, v_{1,2}, v_{2,2}$.

Finally suppose that there is an edge $e$ in $G$  with both of its end vertices in $V(G) \setminus \{v_{1,1}, v_{2,1}, v_{1,2}, v_{2,2}\}$. Therefore the vertices $v_{1,1}, v_{2,1}, v_{1,2}, v_{2, 2}$ are bipartite dominating vertices in $G-e$  and hence $\gmb'(G-e) =2$, a contradiction with $G$ being $2$-$\gmb'$-critical. Therefore at least one end points of each edge of $G$ is a bipartite dominating vertex and thus $\deg_G(x)=2$ for any $x \in V(G) \setminus \{v_{1,1}, v_{2,1}, v_{1,2}, v_{2,2}\}$.
\qed

If $n,m\ge 3$, then let $B_{n,m}$ be the bipartite graph with a bipartition $V_1, V_2$, where $|V_1| = n$, $|V_2| = m$. The set $V_1$ contains exactly two vertices of degree $m$, the set $V_2$ contains exactly two vertices of degree $n$, while all the other vertices of $B_{n,m}$ are of degree $2$. This uniquely defines $B_{n,m}$. With this definition in hand the following characterization of connected bipartite 2-$\gmb'$-critical graphs can be deduced from  Lemmas~\ref{lem:mbd 7} and~\ref{lem:Bipartite2}.

\begin{theorem}
A connected bipartite graph $G$ is $2$-$\gmb'$-critical if and only if either $G = K_{2,m}$, $m\ge 3$, or $G=B_{n,m}$, $n,m\ge 3$.
\end{theorem}

A bit surprisingly, there exist also connected, non bipartite 2-$\gmb'$-critical graphs. Let $G$ be a graph obtained from two disjoint copies of $K_3$ and a vertex $x$ by connecting $x$ with exactly one vertex in each copy of $K_3$. The graph $G$ is 2-$\gmb'$-critical graphs. In view of this example and of Proposition~\ref{Prop:necessary-conditions} and Lemma~\ref{lem:mbd 7}, we conclude the section with: 

\begin{problem}
Characterize $2$-$\gmb'$-critical graphs.  
\end{problem}

\section{SMBD game critical graphs}
\label{sec: SMBD-critical}

In this section, we introduce SMBD game critical graphs. By Lemma~\ref{lem:doky}~(iii), (iv), the SMBD number of a graph never decreases by adding an edge to the graph. Hence, in this case we define critical graphs as follows.

\begin{definition}\label{definition 2} 
If $G$ is a graph and $\tau \in \{\gsmb,\gsmb'\}$, then $G$ is $\tau$-critical, if $\tau(G) < \tau(G+e)$, for any edge $e \notin E(G)$. 
\end{definition}

Note that in view of Lemma~\ref{lem:doky}~(iii) and (iv), if $G$ is $\tau$-critical, where $\tau \in \{\gsmb,\gsmb'\}$, then $\tau(G) < \infty$.  If $G$ is $\tau$-critical and $\tau(G)=k$, where $\tau \in \{\gsmb,\gsmb'\}$, then we say that $G$ is a {\em $k$-$\tau$-critical}. 

Let $\gsmb(G)=k$. Iteratively adding edges $e \notin E(G)$ with the property $\gsmb(G+e)=k$ we arrive at a $k$-$\gsmb$-critical graph. Hence, any graph $G$ with $\gsmb(G)=k$ is a spanning subgraph of a $k$-$\gsmb$-critical graph.  Analogous conclusion holds for a $k$-$\gsmb'$-critical graph. So such critical graphs clearly exists. 

\begin{proposition}\label{prop:1smb}
If $G$ is a graph, then the following holds. 
\begin{enumerate}
    \item[(i)] $G$ is $1$-$\gsmb$-critical if and only if $G = 2K_1$ or $G = K_n\cup 2K_1$, $n \geq 1$.
    \item[(ii)] $G$ is $1$-$\gsmb'$-critical if and only if $G = K_n\cup K_1$, $n \geq 1$.
\end{enumerate}
\end{proposition}

\proof
 (i) Assume that  $G = 2K_1$ or $G = K_n\cup 2K_1$, $n \geq 1$. Since $G$ has at least two isolated vertices, $\gsmb(G)=1$. Moreover, any edge $e\notin E(G)$ is incident to at least one of the isolated vertices. Therefore, $G+e$ contains at most one isolated vertex and hence $\gsmb(G+e)\geq 2$ for any $e\notin  E(G)$. Thus $G$ is $1$-$\gsmb$-critical.
 
 Conversely, assume that $G$ is $1$-$\gsmb$-critical, that is, $\gsmb(G)=1$ and $\gsmb(G+e)\geq 2$ for any $e\notin E(G)$. Hence
 $G$  has at least two isolated vertices, and for any $e \notin E(G)$, the graph $G+e$ has at most one isolated vertex. Therefore, if $e\notin  E(G)$, then $e$ is incident to at least one of the isolated vertices. This is possible only when $G$ is either $2K_1$ or $K_n\cup 2K_1$ for $n \geq 1$.

 \medskip
 (ii) Assume that $G$ is $K_n\cup K_1$ for $n \geq 1$. Since $G$ has an  isolated vertex, $\gsmb'(G)=1$. Moreover, any edge $e\notin E(G)$ is incident to the vertex of $K_1$. Therefore, $G+e$ has no isolated vertices and hence $\gsmb'(G+e)\geq 2$ for any $e\notin  E(G)$. Thus $G$ is $1$-$\gsmb'$-critical.

 Conversely, assume that $G$ is $1$-$\gsmb'$-critical, that is, $\gsmb'(G)=1$ and $\gsmb'(G+e)\geq 2$ for any $e\notin E(G)$. Thus $G$  has at least one isolated vertex and for any $e \notin E(G)$ it follows that $G+e$ has no isolated vertices. Therefore, $G$ has an isolated vertex $x$ such that any edge $e\notin E(G)$ is incident to $x$. This is possible only when $G = K_n\cup K_1$, $n \geq 1$. 
 \qed

To describe $2$-$\gsmb'$-critical graphs and $2$-$\gsmb'$-critical graphs, we define the following two families of graphs. A graph $G_n'$ is a graph obtained from $K_n$ by attaching exactly two pendent vertices to a vertex of $K_n$, and a graph $G_n''$ is obtained from $K_n$ by respectively attaching exactly two pendent vertices to two different vertices of $K_n$. See Fig.~\ref{fig:ilustrative-example}.

\begin{figure}[ht!]
    \begin{center}
        \begin{tikzpicture}[scale=0.6,style=thick,x=1cm,y=1cm]
            \def\vr{3pt}
            \begin{scope}[xshift=0cm, yshift=0cm] 
                \coordinate(g) at (2,0);
                \coordinate(h1) at (2,2);
                \coordinate(h2) at (6,0);
                \coordinate(h3) at (6,2);
                \coordinate(h4) at (4,4);
                \coordinate (h5) at (8,1);
                \coordinate (h6) at (8,3);
                 \coordinate (h5) at (8,1);
                \coordinate (h6) at (8,3);
                 \coordinate (h7) at (1,1);
               	
                \draw (g) -- (h1);
                \draw (g) -- (h2); 
                \draw (g) -- (h3);
                \draw (g) -- (h4);
                \draw (h1) -- (h2);
                \draw (h1) -- (h3); 
                \draw (h1) -- (h4); 
                \draw (h2) -- (h3);
                \draw (h2) -- (h4);
                \draw (h3) -- (h4);
                \draw (h3) -- (h5);
                \draw (h3) -- (h6);
                \draw(g)[fill=white] circle(\vr);
                \foreach \i in {1,...,7}
                { \draw(h\i)[fill=white] circle(\vr);
                }
            
            \end{scope}
            \begin{scope}[xshift=10cm, yshift=0cm] 
                \coordinate(g) at (2,0);
                \coordinate(h1) at (2,2);
                \coordinate(h2) at (6,0);
                \coordinate(h3) at (6,2);
                \coordinate(h4) at (4,4);
                \coordinate (h5) at (8,1);
                \coordinate (h6) at (8,3);
                 \coordinate (h5) at (8,1);
                \coordinate (h6) at (8,3);
                 \coordinate (h7) at (1,1);
                \coordinate (h8) at (1,3);
                \draw (g) -- (h1);
                \draw (g) -- (h2); 
                \draw (g) -- (h3);
                \draw (g) -- (h4);
                \draw (h1) -- (h2);
                \draw (h1) -- (h3); 
                \draw (h1) -- (h4); 
                \draw (h2) -- (h3);
                \draw (h2) -- (h4);
                \draw (h3) -- (h4);
                \draw (h3) -- (h5);
                \draw (h3) -- (h6);
                \draw (h1) -- (h7);
                \draw (h1) -- (h8);
                \draw(g)[fill=white] circle(\vr);
                \foreach \i in {1,...,8}
                { \draw(h\i)[fill=white] circle(\vr);
                }
            \end{scope}

        \end{tikzpicture}
        \caption{Graphs $G_5' \cup K_1$ and $G_5''$}
        \label{fig:ilustrative-example}
    \end{center}
\end{figure}

For the proof of the next result we recall from~\cite[Proposition 4.4]{divakaran-2024} that if $G$ is a graph with $\delta(G) \ge 1$, then $\gsmb(G) = 2$ if and only if $G$ has at least two strong support vertices, and that $\gsmb'(G) = 2$ if and only if $G$ has a strong support vertex.

\begin{theorem}
\label{thm:2-critical-smb}
If $G$ is a graph, then the following holds. 
\begin{enumerate}
    \item[(i)] $G$ is $2$-$\gsmb$-critical if and only if $G = G_n'\cup K_1$, $n \geq 1$, or  $G = G_n''$, $n \geq 2$.
    \item[(ii)] $G$ is $2$-$\gsmb'$-critical if and only if $G = G_n'$, $n\ge 1$.
\end{enumerate}
\end{theorem}

\proof
(i) Assume that $G$ is $2$-$\gsmb$-critical.
That is, $\gsmb(G)=2$ and $\gsmb(G+e) > 2$ for any edge $e\notin E(G)$. Hence $G$ either has exactly one strong support vertex and exactly one isolated vertex, or $\delta(G) \geq 1$ and $G$ has exactly two strong support vertices. 

First, assume that $G$ has exactly one strong support vertex $u$ with two leaf neighbors $u_1,u_2$ and isolated vertex $v$. Since $\gsmb(G+e) > 2$ for any edge $e\notin E(G)$, it follows that $G+e$ has either no strong support vertices or no isolated vertices. Hence any edge $e \notin E(G)$ is either incident with $v$ or with a vertex from $\{u_1,u_2\}$.  Thus $G = G'\cup K_1$. Assume now that $G$ has exactly two strong support vertices and $\delta(G)\geq 1$. Since $G$ is critical, $G+e$ has at most one strong support vertex. This is possible only when any edge $e\notin E(G)$ is incident with at least one pendant vertex of $G$ that eliminates at least one strong support vertex. Thus $G = G_n''$.

Conversely assume that $G = G_n'\cup K_1$, $n \geq 1$, or  $G = G_n''$, $n \geq 2$. Then $G$ has either a strong support vertex and an isolated vertex, or two strong support vertices. Therefore $\gsmb(G)=2$. Since any edge $e\notin E(G)$ is incident to either the isolated vertex or a pendant vertex, the graph $G+e$ has either  at most one strong support vertex and no isolated vertices, or one isolated vertex and no strong support vertices. Therefore $\gsmb(G+e)>2$. Thus $G$ is $2$-$\gsmb$-critical.

(ii) Assume that $G$ is $2$-$\gsmb'$-critical, that is, $\gsmb'(G)=2$ and $\gsmb'(G+e) > 2$ for any edge $e\notin E(G)$. Hence $\delta(G)=1$ and $G$ has exactly one strong support vertex. Moreover, any edge $e\notin E(G)$ is incident to a pendant vertex so that $G+e$ has no strong support vertices. Hence $G = G_n'$.

Conversely assume that $G = G_n'$. Since $G$ has a strong support vertex with exactly two leaf neighbors, $\gsmb'(G)=2$. As any edge $e\notin E(G)$ is incident with at least one of the pendant vertices, $G+e$ has no strong support vertices and hence $\gsmb'(G+e) > 2$. Thus $G$ is $2$-$\gsmb'$-critical. 
\qed

Recall that $G|X$ denotes the graph $G$ in which vertices from $X\subseteq V(G)$ are considered as being already dominated. For $v \in V(G)$, let $G_v$ be a graph obtained from $G|N[v]$ by deleting all dominated vertices $x \in V(G|N[v])$ with $N[x] \subseteq N[v]$.  Assume now that $G$ is 2-$\gsmb$-critical, so that in view of Theorem~\ref{thm:2-critical-smb}, $G$ is $K_1 \cup G_n'$ for some $n \geq 1$, or $G_n''$ for $n \geq 2$. Then there exists a vertex $v \in V(G)$ (an isolated vertex or a strong support vertex of $G$) such that $G_v$ is 2-$\gsmb'$-critical. We wonder if this holds also for $k$-$\gsmb$-critical graphs with $k \geq 3$:

\begin{question}
Let $G$ be a $k$-$\gsmb$-critical graph, $k\ge 3$. Is it true that there exists a vertex $v \in V(G)$ such that $G_v$ is $\gsmb'$-critical?    
\end{question}

\section*{Acknowledgements}

T.\ Dravec and S.\ Klav\v{z}ar have been supported by the Slovenian Research Agency ARIS (research core funding P1-0297 and projects N1-0285, N1-0355). A lot of work on this paper has been done during the Workshop on Games on Graphs II, June 2024, Rogla, Slovenia, the authors thank the Institute of Mathematics, Physics and Mechanics, Ljubljana, Slovenia for supporting the workshop.  

\section*{Declaration of interests}
 
The authors declare that they have no conflict of interest. 

\section*{Data availability}
 
Our manuscript has no associated data.

\end{document}